\documentclass[brochure,12pt,english]{bourbaki}
\usepackage[matrix,arrow]{xy}
\newtheorem{hypo}[defi]{Hypothesis}
\newcommand{\qqed}{\hspace*{\fill}$\Box$}
\newcommand{\End}{{\rm End}}
\newcommand{\IC}{{\mathbb C}}
\newcommand{\IP}{{\mathbb P}}
\newcommand{\IQ}{{\mathbb Q}}
\newcommand{\IR}{{\mathbb R}}
\newcommand{\IZ}{{\mathbb Z}}
\date{Novembre 2005}
\bbkannee{58\`eme ann\'ee, 2005-2006}
\bbknumero{954}
\title{PROJECTIVITY OF K\"AHLER MANIFOLDS\\ - KODAIRA'S PROBLEM
}
\subtitle{after C.\ Voisin}
\author{Daniel HUYBRECHTS}
\address{Universit\"at Bonn\\
Mathematisches Institut\\
Beringstra\ss e 1\\
D--53115 Bonn, Germany}
\email{huybrech@math.uni-bonn.de}
\begin{document}
\maketitle


\bigskip
There are various geometric structures that can be studied on a
topological manifold $M$. Depending on one's geometric taste, it
is important to know whether $M$ can be endowed with  a symplectic
form, whether (special) Riemannian metrics can be found or
whether $M$ carries an algebraic structure. Often, the existence
of a certain geometric structure imposes topological conditions
on $M$. In other words, it may happen that a given topological
manifold does simply not allow one's favorite geometry. E.g.\ if
$M$ is compact and $b_2(M)=0$ the manifold $M$ cannot be
symplectic, or if $b_1(M)=1$ no K\"ahler metrics can exist.

In order to fully understand the relation between two sorts of
geometries, it is important to know whether they impose the same
topological obstructions. In other words, does the existence of
one of the two on a given manifold topological $M$ imply the
existence of the other one? This is a report on the work of
Claire Voisin \cite{V1,V2} that sheds light on an old question,
usually attributed to Kodaira, that asks for the topological
relation between K\"ahler geometry and projective geometry.

\smallskip

In the following we let $M$ be a compact manifold that can be
endowed with the structure of a complex manifold. Once a complex
structure is chosen, one studies Riemannian metrics $g$ that are
`compatible' with it. One possible compatibility condition is to
require that $g$ be hermitian, i.e.\ that the complex structure
thought of as an almost complex structure $I$ is orthogonal with
respect to $g$. It is not difficult to see that a hermitian
structure can always be found. It is, however, a completely
different matter to find a hermitian structure $g$ such that its
fundamental form $\omega:=g(I~~,~~)$ is closed, i.e.\ $g$
satisfies the K\"ahler condition. Indeed, the classical theory of
K\"ahler manifolds shows that the existence of a K\"ahler metric
imposes strong conditions on the topology of $M$, which are not
satisfied by arbitrary complex or symplectic manifolds. For
instance, the odd Betti numbers of a compact K\"ahler manifold
are even, K\"ahler manifolds are formal and  their fundamental
groups satisfy further conditions. (In contrast, if only one of
the two structures, complex or symplectic, is required, then any
finitely presentable group can be realized.)

On the other hand,  K\"ahler manifolds are quite common. Indeed,
any complex submanifold of the complex projective space $\IP^n$
admits a K\"ahler metric - the restriction of the Fubini--Study
metric is an example.  Conversely, one might wonder whether a
compact complex manifold that admits a K\"ahler structure can
always be realized as a complex submanifold of $\IP^n$ or, in
other words, whether the complex structure is projective. This is
obviously not the case, general complex tori $\IC^n/\Gamma$
($n\geq2$)  and general K3 surfaces provide counter-examples. In
fact, a famous theorem of Kodaira proves that a K\"ahler manifold
is projective if and only if the K\"ahler metric can be chosen
such that the cohomology class of its fundamental form $\omega$
is integral, i.e.\ $[\omega]\in H^2(X,\IZ)$ (see
\cite[Thm.4]{Kod1}).

In these examples one observes that although the given complex
structure is not projective, it becomes projective after a small
deformation. Kodaira proved that in fact any K\"ahler surface can
be deformed to a projective surface (see \cite[Thm.23]{Kod0} and
\cite{Kod2}). Thus, as deforming the complex structure does not
change the diffeomorphism type of the manifold, there is no
topological difference between compact K\"ahler surfaces and
algebraic surfaces. (Let us also mention that in fact any compact
surface $X$ with even $b_1(X)$ is K\"ahler, i.e.\ for surfaces
the condition to be K\"ahler is a topological condition. This
fails in higher dimensions, due to a famous example of Hironaka
\cite{Hir} of a compact K\"ahler manifold that deforms to complex
manifold which is no longer K\"ahler.) Note in passing that a
similar results holds true for symplectic manifolds: Clearly, any
given symplectic form $\omega$ can be deformed to a symplectic
form with integral cohomology class.

Kodaira's problem, which apparently has never been stated by
himself in this form,  asks for the higher-dimensional version of
his result: {\it Can any compact K\"ahler manifold be deformed to
a projective manifold?}

More in the spirit of the general philosophy explained above, one
could ask whether the topological manifold underlying a compact
K\"ahler manifold may also be endowed with the structure of a
projective manifold. This question had been open for a very long
time.  As Kodaira's arguments to prove the two-dimensional case
use a great deal of classification theory of surfaces, there was
little hope to generalize them to higher dimensions.


Recent work of Claire Voisin fills this gap \cite{V1,V2,V3}. She
succeeded in showing that topology makes a difference between
compact K\"ahler manifolds and those that are projective. In
other words, there exist compact topological manifolds that admit
the structure of a K\"ahler manifold without carrying also the
structure of a projective manifold. More precisely, Voisin shows
the stronger statement:

\begin{theo}[\cite{V1}]\label{VoisinMainThm}
In any dimension $\geq4$ there exists a compact K\"ahler manifold
$X$ whose rational cohomology ring $H^*(X,\IQ)$ cannot be realized
as the rational cohomology ring of a projective manifold.
\end{theo}

Voisin originally worked with the integral cohomology ring
$H^*(X,\IZ)$, but Deligne then pointed out the stronger version
above.


%
One could wonder whether the answer to these questions would be
different if the topological manifold satisfies further
conditions, e.g.\ if it is in addition simply-connected. Some of
these questions have been addressed and answered by Voisin in
\cite{V1,V2} and we will comment on them on the way.

\bigskip

Although the examples are obtained by particular constructions, the principal
ideas of \cite{V1,V2} are of a more general nature and might be applicable
in other situations.

The $i$-th cohomology of a compact K\"ahler manifold is naturally
endowed with a Hodge structure of weight $i$, which can be
polarized (on the primitive part) if the manifold is projective.
The idea is to show that there exist compact K\"ahler manifolds
whose cohomology does not  admit Hodge structures that are
compatible with both, the given cup-product and a polarization
Roughly,
 there are three steps {\bf A-C}, the first
two of which are purely Hodge-theoretical and only the last one
has a geometric flavor.

\smallskip

{\bf (A)} Certain algebraic structures on a rational vector space
$A$ are not compatible with any {\it polarizable} Hodge structure
(of weight $k$) on $A$.

\begin{rema}
In the examples, the algebraic structure will  be a specific
endomorphism $\Phi:A\to A$, but others are in principle possible.
That the algebraic structure is not compatible with any
polarizable Hodge structure means in the case of an endomorphism
$ \Phi$ that one cannot find a Hodge structure on $A$ such that
$\Phi$ becomes an endomorphism of it and such that the Hodge
structure can be polarized. 
\end{rema}

{\bf(B)} Suppose $\bigoplus H^\ell$ is a graded $\IQ$-algebra
whose direct summands $H^\ell$ are Hodge structures of weight
$\ell$ and such that the multiplications $H^{\ell_1}\otimes
H^{\ell_2}\to H^{\ell_1+\ell_2}$ are homomorphisms of Hodge
structures. Suppose furthermore that this $\IQ$-algebra structure
allows us to detect a subspace $A\subset H^k$ such that: i)
$A\subset H^k$ is a Hodge substructure.
ii) An algebraic structure as in {\bf (A)} is compatible with this Hodge structure.\\
Then $H^k$ does not admit a polarization.

\begin{rema}
Subspaces that are defined purely in terms of the $\IQ$-algebra
structure do define Hodge substructures. We shall also need a
refined version of this, which is due to Deligne.

The compatibility in ii) is more difficult to check, but relies
on the same principle. For an endomorphism $\Phi$ the idea goes
as follows: Firstly, find two Hodge substructures $A,A'\subset
H^k$ and a Hodge substructure $\Delta\subset A\oplus A'\subset
H^k$ which is the graph of an isomorphism $A\cong A'$. Secondly,
prove that under the induced isomorphism of Hodge structures  $
A\oplus A\cong A\oplus A'$ the graph of $\Phi$ is a Hodge
substructure.
\end{rema}

{\bf (C)} Construct compact K\"ahler manifolds
such that the above principles apply
to its cohomology ring $\bigoplus H^\ell(X,\IQ)$.
Then $H^*(X,\IQ)$ should not be realizable by a smooth projective variety.

\begin{rema}
This works best for Hodge structures of weight one ($k=1$). In
this case $H^1(X,\IQ)$ of a smooth projective variety $X$ admits a
polarized Hodge structure. For the Hodge structure of weight two
on $H^2(X,\IQ)$ one needs an extra argument, for only the
primitive part of it admits a polarization.
\end{rema}

\smallskip

This report roughly follows these three steps. Some of the
algebraic structures in Section \ref{ImpSect} might seem rather
ad hoc, as their geometric origin is only explained in Section
\ref{Exassect}. However, I found it helpful for my own
understanding to completely separate the arguments that explain
why certain $\IQ$-algebras cannot be realized as the cohomology of
a projective manifold from the part that  contains the
construction of compact K\"ahler manifolds that do realize these
$\IQ$-algebras.

\smallskip

{\bf Acknowledgements:} I wish to thank Claire Voisin for
patiently answering my questions and for her valuable comments on
a first draft of these notes. I am grateful to C.-F.\
B\"odigheimer, U.\ G\"ortz, M.\ Lehn, P.\ Stellari,
J.\ Stix, R.\ Thomas, B.\ Totaro, and T.\ Wedhorn for their help, comments,
and suggestions.

\section{Hodge structures (of weight one and two)}
\subsection{Recollections}

A \emph{Hodge structure of weight} $k$ on a $\IQ$-vector space
$A$ is given by a direct sum decomposition
\begin{equation}\label{HS}
A_\IC:=A\otimes_\IQ\IC={\displaystyle\bigoplus_{p+q=k}A^{p,q}}~~{\rm
such~that}~~\overline{A^{p,q}}=A^{q,p}.\end{equation}

A direct sum decomposition (\ref{HS}) can also be described in
terms of a representation $\rho:\IC^*\to{\rm Gl}(A_\IR)$ such
that the $\IC$-linear extension of $\rho(z)$ satisfies
$\rho(z)|_{A^{p,q}}=z^p\bar z^q\cdot{\rm id}$. The \emph{Hodge
classes} of a Hodge structure of weight $2k$ on $A$ are the
elements in $A^{k,k}\cap A$.

We shall be particularly interested in Hodge structures of weight
one and two.

\begin{rema} Recall that Hodge structures of weight one with $A^{p,q}=0$ for
$pq\ne0$ which are integral, i.e.\  $A=\Gamma_\IQ$ for some
lattice $\Gamma$, are in bijection with complex tori. Indeed, to
a Hodge structure of weight one on $\Gamma_\IQ$ given by
$\Gamma_\IC=A^{1,0}\oplus A^{0,1}$ one associates the complex
torus $A^{1,0}/\Gamma$, where $\Gamma$ is identified with its
image under the projection $A_\IC\to A^{1,0}$.
\end{rema}

A $\IQ$-linear map $\varphi:A\to A'$ is a \emph{morphism (of
weight $m$) of Hodge structures}
 $$A_\IC={\displaystyle\bigoplus_{p+q=k}}A^{p,q}~~{\rm
and}~~A'_\IC={\displaystyle\bigoplus_{r+s=\ell}}A'^{r,s}$$ of
weight $k$ and $\ell=k+2m$, respectively, if
$\varphi(A^{p,q})\subset A'^{p+m,q+m}$. If the two Hodge
structures correspond to $\rho:\IC^*\to{\rm Gl}(A_\IR)$ and
$\rho':\IC^*\to{\rm Gl}(A'_\IR)$, respectively, then this
condition is equivalently expressed by $\varphi(\rho(z)v)=
|z|^{2m}\rho'(z)\varphi(v)$ for all $v\in A$ and $z\in\IC^*$.

A \emph{Hodge substructure} of a Hodge structure of weight $k$ on
$A$ is given by a subspace $A'\subset A$ such that
$A'_\IC=\bigoplus \left(A^{p,q}\cap A'_\IC\right)$ or,
equivalently, such that $A'_\IC\subset A_\IC$ is in\-variant under
the representation $\rho:\IC^*\to{\rm Gl}(A_\IR)$ that
corresponds to the given Hodge structure on $A$.

The \emph{tensor product} $A\otimes_\IQ A'$ of two $\IQ$-vector
spaces $A$ and $A'$ endowed with Hodge structures of weight $k$
and $\ell$, respectively, comes with a natural Hodge structure of
weight $(k+\ell)$:
$$(A\otimes_\IQ
A')^{r,s}:={\displaystyle\bigoplus_{p+p'=r,q+q'=s}}
A^{p,q}\otimes_\IC A'^{p',q'}.$$ In other words, the Hodge
structure is given by $\rho\otimes\rho'$.

Note that $A_2:=\bigwedge^2A_1$ of a Hodge structure of weight one
$A_1$ is naturally a Hodge structure of weight two with
$A_2^{2,0}:=\bigwedge^2A_1^{1,0}$, $A_2^{1,1}:=A_1^{1,0}\otimes
A_1^{0,1}$, and $A_2^{0,2}:=\bigwedge^2A_1^{0,1}$.

A \emph{polarization} of a Hodge structure of weight one
$A_\IC=A^{1,0}\oplus A^{0,1}$ is a skew-symmetric form $q\in
\bigwedge^2 A^*$ such that
\begin{equation}\label{HR}\xymatrix{A_\IC\times
A_\IC\ar[r]& \IC,}~~\xymatrix{(v,w)\ar@{|->}[r]&iq(v,\overline w)}
\end{equation}
(where $q$ is extended $\IC$-linearly) satisfies the
Hodge--Riemann relations:\\
i) $A^{1,0}$ and $A^{0,1}$ are orthogonal with respect to
(\ref{HR}).\\
ii) The restriction of (\ref{HR}) to  $A^{1,0}$ and to $A^{0,1}$
is positive, respectively negative, definite.

\begin{rema}
With this definition a polarization is always rational.
Furthermore, the form $q$ considered as an element of the induced
weight-two Hodge structure on $\bigwedge^2A^*$ is of type $(1,1)$.
Since it is rational, $q$ is a Hodge class (of weight two). Note
that any Hodge substructure of a weight-one polarized Hodge
structure is naturally polarized.
\end{rema}

\begin{exem}
Let $X$ be a compact K\"ahler manifold of dimension $n$. The
Hodge decomposition $$H^1(X,\IC)=H^{1,0}(X)\oplus H^{0,1}(X)$$
defines a Hodge structure of weight one on $H^1(X,\IQ)$.

Suppose $X$ is projective and $\omega\in H^2(X,\IZ)$ is the class of a
hyperplane section, then $q(\alpha)=\int_X\alpha^2\omega^{n-1}$
is a polarization of the natural Hodge structure of weight one on
$H^1(X,\IQ)$.

If we drop the condition that $q$ be rational, then any K\"ahler
class on a compact K\"ahler manifold $X$ would yield a form on
the Hodge structure of weight one on $H^1(X,\IQ)$ that satisfies
the Hodge--Riemann relations i) and ii).
\end{exem}

The notion of a polarization exists for Hodge stuctures of
arbitrary weight, but we shall only need it for weight
one, explained above, and for weight two. For a Hodge structure
of weight two $A_\IC=A^{2,0}\oplus A^{1,1}\oplus A^{0,2}$ a
\emph{polarization} is a symmetric bilinear form $q\in S^2A^*$ such that:\\
i) The $A^{p,q}$ are pairwise orthogonal with respect to $(v,w)\mapsto
q(v,\overline w)$.\\
ii) For $0\ne v\in A^{p,q}$ one has $-i^{p-q}q(v,\overline v)>0$.

\begin{exem}\label{exaHSwttwo}
If $X$ is compact K\"ahler of dimension $n$, then
$H^2(X,\IQ)$ comes with a natural Hodge structure of weight two
$H^2(X,\IC)=H^{2,0}(X)\oplus H^{1,1}(X)\oplus H^{0,2}(X)$ given
by the Hodge decomposition. If $X$ is projective and $\omega\in
H^2(X,\IZ)$ is the class of a hyperplane section, then
$$q(\alpha)=\int_X\alpha^2\omega^{n-2}$$defines a polarization on
the \emph{primitive cohomology}
$$H^2(X,\IQ)_{\rm p}:=\{\alpha\in
H^2(X,\IQ)~|~\alpha\wedge\omega^{n-1}=0\}.$$

Note that due to the Hodge--Riemann bilinear relation
$H^{1,1}(X,\IR)\cong H^{1,1}(X,\IR)_{\rm p}\oplus\IR\omega$ does
not contain any $q$-isotropic subspace of dimension $\geq2$. Also,
$H^2(X,\IR)$ does not contain Hodge substructures of dimension
$\geq2$ which are $q$-isotropic.
\end{exem}

\subsection{Detecting Hodge structures
algebraically}\label{DelSect}

The following observation is the key to a general principle, due
to Deligne, which allows one to identify Hodge substructures
algebraically.

\begin{lemm}\label{DeligneLem}
 Let $H_\IC=\bigoplus_{p+q=k} H^{p,q}$ be a Hodge
structure of weight $k$ on a $\IQ$-vector space $H$ given by a
representation $\rho:\IC^*\to{\rm Gl}(H_\IR)$ and let $Z\subset
H_\IC$ be an algebraic subset which is invariant under
$\rho(\IC^*)$. Suppose the span $\langle Z'\rangle$ of an
irreducible component $Z'\subset Z$  is of the form
$H'\otimes_\IQ\IC$ with $H'\subset H$ a $\IQ$-subspace. Then $H'$
is a Hodge substructure of $H$.
\end{lemm}

\begin{proof}
Since $\IC^*$ is connected, the $\IC^*$-action leaves invariant
the irreducible components of $Z$. Hence, also $\langle
Z'\rangle$ is $\IC^*$-invariant. For  $\langle
Z'\rangle=H'\otimes_\IQ\IC$ this is equivalent to saying that
$H'\subset H$ is a Hodge substructure.
\end{proof}

In \cite{V1,V2} the lemma is applied in various situations. The
algebraic set $Z$ is always defined by algebraic conditions on
homomorphisms of Hodge structures and thus automatically
invariant under $\IC^*$. Usually, one starts with several Hodge
structures of weight $\ell$ on $\IQ$-vector spaces $H^\ell$ and
homomorphisms of Hodge structures $H^{\ell_1}\otimes H^{\ell_2}\to
H^{\ell_1+\ell_2}$, $a\otimes b\mapsto a\cdot b$. (Think of the
cohomology of a smooth projective variety  or of a compact
K\"ahler manifold.)

We shall in particular encounter algebraic subsets of the form

$$Z_1=\{\alpha\in H^k_\IC~|~ \alpha^2=0\}~~~{\rm or}~~~
Z_2=\left\{a\in H^k_\IC~|~{\rm
rk}\left(\xymatrix{\!H^{\ell}_{\IC}\ar[r]^-{a\cdot}&H^{k+\ell}_\IC}\!\right)\leq
m\right\}.$$

Let us sketch the argument that shows that these sets are
$\IC^*$-invariant in the example $Z=Z_2$. By definition of the
Hodge structure on $H^{\ell_1}\otimes H^{\ell_2}$ and the
hypothesis that the multiplication $a\otimes b\mapsto a\cdot b$
is a morphism of Hodge structures, one has $\rho(z)(a)\cdot
b=\rho(z)(a\cdot(\rho(z^{-1})(b)))$. Thus, the endomorphism given
by multiplication with $\rho(z)(a)$ and $a$, respectively, differ
by automorphisms $\rho(z)\in{\rm Gl}(H^{k+\ell}_\IR)$ and
$\rho(z^{-1})\in{\rm Gl}(H^{\ell}_{\IR})$. In particular, ${\rm
rk}(\rho(z)(a)\cdot)={\rm rk}(a\cdot )$ and hence $a\in Z$ if and
only if $\rho(z)(a)\in Z$.

Note that it might well happen that $\langle Z\rangle$ is defined
over $\IQ$, but not $\langle Z'\rangle$.


\medskip

 Let us illustrate the use of Deligne's principle in a
concrete situation that will be at the heart of the subsequent
discussion. Suppose we are given  a graded $\IQ$-algebra
$\bigoplus H^k$, an integer $\ell\in\IZ$ and a subspace  $0\ne
H'\subset H^\ell$. Then define for $i\geq1$ the  $\IQ$-subspace
\begin{equation}\label{DFNP}
P_i:=\left\{a\in
H^2~|\left(\!\xymatrix{\bigotimes^{i}H'\ar[r]^-{\cdot a}&H^{\ell
i+2}}\right)\!=0\right\}.
\end{equation}
We shall later fix in addition an integer $m>1$ and consider the
two subspaces $$P_1\subset P_m\subset H^2$$
and
the algebraic subset of $P_{m\IC}$:
\begin{equation}\label{DFNZ}
Z:=\{a\in P_{m\IC}~|~{\rm
Ker}\left(\!\xymatrix{H'_\IC\ar[r]^{\cdot
a}&H^{\ell+2}_\IC}\right)\ne0\}. \end{equation} Then $Z$ contains
$P_{1\IC}$ and we denote its image in $(P_m/P_1)_\IC$ by $\bar Z$
(which is again algebraic). Furthermore, let $e\in Z\cap P_m$ be
such that $\IC\bar e\subset \bar Z$ is an irreducible component
of $\bar Z$.

\begin{coro}\label{DelCor} Suppose  each $H^k$ is endowed with a
Hodge structure of weight $k$ such that the multiplications are
morphisms  of Hodge structures and such that $H'\subset H^\ell$
is a Hodge substructure. Then

 {\rm i)} the $P_i\subset H^2$ are Hodge
substructures,

{\rm ii)} the element $\bar e\in P_m/P_1$ is of type $(1,1)$,
i.e.\ a Hodge class, and

{\rm iii)} ${\rm Ker}(H'\stackrel{\cdot e}{\longrightarrow}
H^{\ell+2})$ is a Hodge substructure of $H^\ell$.
\end{coro}

\begin{proof} The $P_i$ can be viewed as the kernels
of the morphisms of Hodge structures $H^2\to
\left(\bigotimes^iH'\right)^*\otimes H^{\ell i+2}$ and are,
therefore, Hodge substructures of $H^2$.

Deligne's principle shows that $\IQ\bar e\subset P_m/P_1$ is a
Hodge substructure. Since any weight two Hodge structure of rank
one is of pure type, one finds $\bar e\in (P_m/P_1)^{1,1}$.

In order to prove iii), use the morphism of Hodge structures
$P_m/P_1\otimes H'\to H^{\ell+2}$.
\end{proof}

\begin{rema}
i) The actual description of $P_m$ is of no importance here. We
only used $P_1\subset P_m$ and the condition on $e$. Note that
$e\in P_m$ itself might be of mixed type, e.g.\ it could be
arbitrarily modified by rational classes in $P_1^{2,0}\oplus
P_1^{0,2}$.

ii) In the applications only the cases $\ell=1$ and $\ell=2$ will
be considered and, moreover,   for $\ell=1$ we will have $H'=H^1$.
\end{rema}


%
%
%
%
%

\section{The impossible ones}\label{ImpSect}
The aim is to exhibit two specific Hodge structures of weight one
respectively two which resist polari\-zation. Section
\ref{ImpSect1} explains Step {\bf A} of the program, whereas
Section \ref{ImpSect2} corresponds to Step {\bf B}.

\subsection{Special endomorphisms excluding polarization}\label{ImpSect1}

Let us start out with an endomorphism $\Phi\in\End(A)$ of a
$\IQ$-vector space $A$ of dimension $2n$. For any field
$\IQ\subset K$ we shall denote by $\Phi_K$ its $K$-linear
extension. We  also use the naturally induced endomorphisms
$\Phi^*$ and $\bigwedge^2\Phi^*$ of $A^*$ and $\bigwedge^2A^*$
respectively.

Denote the set of all eigenvalues of $\Phi$ by
$EV(\Phi):=\{\mu_1,\ldots,\mu_{2n}\}$ and by $K_\Phi$
the splitting  field of the characteristic
polynomial of $\Phi$, i.e.\ $K_\Phi=\IQ(\mu_1,\ldots,\mu_{2n})$.

Henceforth, we shall assume that:\\
\begin{hypo}\label{hypo}
{\rm i)} $\mu_i\not\in\IR$ for all $i$, and {\rm ii)} $G:={\rm
Gal}(K_\Phi/\IQ)$ acts as the symmetric group $S_{2n}$ on
$EV(\Phi)$.
\end{hypo}

\medskip

\begin{exem} It is not difficult to find explicit examples of
endomorphisms $\Phi$ satisfying these conditions:\\
 -- Let $A=\IQ^2$, hence $n=1$, and
$\Phi=\left(\begin{matrix}0&1\\-1&0\end{matrix}\right)$.
Then $\{\mu_1,\mu_2\}=\{\pm i\}$.\\
-- Let $A=\IQ^4$ and $\Phi=\left(\begin{matrix}0&0&0&-1\\
 1&0&0&1\\
 0&1&0&0\\
 0&0&1&0\end{matrix}\right)$. The characteristic polynomial
of $\Phi$ is $x^4-x+1$ whose Galois group is the symmetric group
(see \cite[Ch.14.6]{Artin}) and which clearly has no real
eigenvalues.
\end{exem}

\begin{rema}\label{Galoispol}
Clearly, ii) implies that $\Phi_\IC\in\End(A_\IC)$ can be
diagonalized. It also yields
$\mu_{i_1}\cdot\ldots\cdot\mu_{i_k}\ne\mu_{j_1}\cdot\ldots\cdot\mu_{j_k}$
for any two distinct multi-indices $i_1<\ldots<i_k$ and
$j_1<\ldots<j_k$.
\end{rema}

\begin{lemm}
Under the assumptions of \ref{hypo} the induced endomorphism
$\bigwedge^k\Phi\in\End(\bigwedge^kA)$ does not admit any
non-trivial invariant subspace.
\end{lemm}

\begin{proof}
Clearly, the eigenvalues of $\bigwedge^k\Phi$ are
$\mu_{i_1}\cdot\ldots\cdot\mu_{i_k}$, $i_1<\ldots <i_k$. Thus, if
$W\subset\bigwedge^kA$ is invariant under $\bigwedge^k\Phi$, then
the eigenvalues of $\psi:=\bigwedge^k\Phi|_W$ are also of this
form. In particular, also $\psi$ can be diagonalized over
$K_\Phi$. Suppose $W\ne0$. Then there exists an eigenvector $v\in
W_{K_\Phi}$ with eigenvalue say $\mu_1\cdot\ldots\cdot \mu_k$.

Being defined over $\IQ$, the
extension of $\psi$ (and of $\bigwedge^k\Phi$) to an endomorphism of
$W_{K_\Phi}$ (respectively $\bigwedge^kA_{K_\Phi}$)
commutes with the action of the
Galois group $G$ on the scalars $K_\Phi$.
Hence, with $\mu_1\cdot\ldots\cdot\mu_k$
also $\mu_{\sigma(1)}\cdot\ldots\cdot\mu_{\sigma(k)}$ is an
eigenvalue of $\psi$ for any $\sigma\in G$.

By Remark \ref{Galoispol}, this shows that all
$\mu_{i_1}\cdot\ldots\cdot\mu_{i_k}$, $i_1<\ldots<i_k$, which are
pairwise distinct, occur as eigenvalues of $\psi$. Hence, $\dim
(W)=\dim(\bigwedge^kA)$ or, equivalently, $W=\bigwedge^kA$.
\end{proof}


\begin{prop}\label{PropnoHodgeclwttwo}
Suppose  $\bigwedge^2\Phi$ respects a  Hodge structure of weight
two  on $A_2:=\bigwedge^2A$ given by
$\bigwedge^2A_\IC=A_2^{2,0}\oplus A_2^{1,1}\oplus A_2^{0,2}$ with
$A_2^{2,0}\ne0$. If $\Phi$ satisfies \ref{hypo}, then
$$A_2^{1,1}\cap {\bigwedge}^2A=\{0\},$$ which is equivalent to saying
that all Hodge classes of $A_2$ are trivial.
\end{prop}

\begin{proof}
As $\bigwedge^2\Phi_\IC$ preserves the bidegree $(p,q)$ of
elements in $\bigwedge^2A_\IC$, the rational subspace
$W:=A_2^{1,1}\cap\bigwedge^2A$ is
$\bigwedge^2\Phi$-invariant. Due to the lemma one either
has $W=\bigwedge^2A$, which is excluded by $A_2^{2,0}\ne0$,
or $W=0$, which proves the assertion.
\end{proof}

\begin{coro}\label{CorwtonenoPo} Suppose $n\geq2$.
A Hodge structure of weight one $A_\IC=A^{1,0}\oplus
A^{0,1}$ that is preserved by $\Phi_\IC$ does not admit a
polarization.
\end{coro}

\begin{proof}
A polarization of the Hodge structure $A_\IC=A^{1,0}\oplus
A^{0,1}$  would be given by a special Hodge class $q$ in the
induced Hodge structure of weight two on
$\bigwedge^2A^*$. However, there are no
non-trivial ones due to the proposition. (Use that
$\Phi^*$ as well satisfies \ref {hypo}.) The assumption $n\geq2$ is needed
in order to ensure that $A_2^{2,0}\ne0$.
\end{proof}

\begin{rema}
Observe that $\Phi$ preserves the Hodge structure if and only if
its graph $\Gamma_\Phi\subset A\oplus A$ is a Hodge substructure.
\end{rema}

\begin{exem}\label{TorusPhi} If $\Phi$ satisfies  i) and ii)
of \ref{hypo}, one
easily constructs Hodge structures of weight one that are
preserved by $\Phi$. This will be needed when it actually comes
to constructing examples.

Pick $n$ distinct eigenvalues $\lambda_1,\ldots,\lambda_n\in
EV(\Phi)$ such that $\lambda_i\ne\bar\lambda_j$ for all $i,j$
(note that due to i) no eigenvalue is real) and let
$A^{1,0}=\bigoplus_{i=1}^n\IC v_i$, where the $v_i\in A_\IC$ are
eigenvectors with eigenvalue $\lambda_i$.

With $\Phi$ being defined over $\IQ$, the complex conjugate
$\bar\lambda$ of an eigenvalue $\lambda\in EV(\Phi)$ is again an
eigenvalue.
Thus, with $A^{0,1}:=\overline{A^{1,0}}$
one has $A_\IC=A^{1,0}\oplus A^{0,1}$.
\end{exem}
\subsection{Identifying the special endomorphisms algebraically}\label{ImpSect2}

We continue the discussion of Section \ref{DelSect} and combine it
with endomorphisms $\Phi$ of the type studied in Section \ref{ImpSect1}.

So, let us consider  a $\IQ$-vector space $A$ of dimension $2n\geq 4$
together with an endomorphism $\Phi$ and let
$H^*=\bigoplus_{k=0}^{4n} H^k$ be a graded $\IQ$-algebra.

To bring both structures together, we assume that there is a
graded inclusion
$${\bigwedge}^*(A\oplus A)\subset H^*$$
satisfying the following conditions. (We shall apply Corollary
\ref{DelCor} with $\ell=1$, $m=4n-2$, and $H'=H^1$.)

\begin{hypo}\label{hypo2}
{\rm i)} $A\oplus A=H^1$,

{\rm ii)} $H^2=\bigwedge^2(A\oplus A)\oplus P\oplus R$, where
$P:=P_{4n-2}$ is defined as in (\ref{DFNP}) and $R$ is some
subspace,

{\rm iii)} $P=P_1\oplus \bigoplus_{i=1}^4e_i\IQ$, and

{\rm iv)} The kernel of the multiplication $H^1\stackrel{\cdot
e_i}{\longrightarrow} H^3$, for $i=1,\ldots,4$,
 equals the subspaces $A\oplus\{0\}$, $\{0\}\oplus A$, $\Delta:=\{
(a,a)~|~a\in A\}$, and the graph $\Gamma_{\Phi}$ of $\Phi$,
respectively. The sum $\sum{\rm Im}(\cdot e_i)\subset H^3$ is
direct.
\end{hypo}

\begin{rema}
Roughly, $e_1$ and $e_2$ will be used to detect certain Hodge
substructures, $e_3$ to identify them, and $e_4$ to view $\Phi$
as a homomorphism between them. The auxiliary space $R$ is later
only needed in order to construct odd-dimensional examples. Due to
Remark \ref{avoidP'Rem} one could even restrict to the case
$P_1=0$.
\end {rema}

\begin{prop}\label{Propnotproj}
Suppose $H^*$ and $\Phi$ meet the conditions of \ref{hypo2} and
\ref{hypo}, respectively. Then $H^*$ cannot be realized as the
rational cohomology ring $H^*(X,\IQ)$ of a projective manifold
$X$.
\end{prop}

\begin{proof}
Suppose  $X$ is a projective manifold that does realize $H^*$. In
the following we will simply identify $H^*(X,\IQ)$ with $H^*$.
Thus, each $H^k$ inherits the  natural Hodge structure of weight
$k$ from $H^k(X,\IQ)$ and the multiplications $H^{\ell_1}\otimes
H^{\ell_2}\to H^{\ell_1+\ell_2}$ are morphisms of Hodge
structures.

Corollary \ref{DelCor} applies and shows that $A\oplus\{0\}$,
$\{0\}\oplus A$, $\Delta$, and the graph $\Gamma_{\Phi}$ are
Hodge substructures of $H^1(X,\IQ)$. Indeed, the only thing that
needs to be checked is that the $\IC \bar e_i$ define irreducible
components of $\bar Z\subset (P/P_1)_\IC$ (the image of $Z$ as in
(\ref{DFNZ})). This follows  from iv): Suppose $\sum a_ie_i\in
Z$. Then there exists $0\ne a\in H^1$ that is annihilated by it.
Thus, $a_i(a\cdot e_i)=0$ for $i=1,\ldots,4$. If e.g.\ $a_i\ne0\ne
a_j$, then $a\in {\rm Ker}(\cdot e_i)\cap{\rm Ker}(\cdot e_j)$.
The description of the kernels shows that this is impossible.

With the identification  of the two Hodge structures on
$A\oplus\{0\}$ and $\{0\}\oplus A$ via $\Delta$, the graph
$\Gamma_{\Phi}$ allows to view $\Phi$ as an endomorphism of the
Hodge structure on $A\oplus\{0\}$.

By Corollary \ref{CorwtonenoPo} this Hodge structure does not
admit a polarization. Hence, also the Hodge structure
$H^1(X,\IQ)$, of which $A\oplus\{0\}$ is a Hodge substructure,
cannot be polarized. This yields a contradiction to the
projectivity of $X$.
\end{proof}

\bigskip

We shall next present a similar result based on an analysis
of Hodge structures of weight two.

Let as before $A$ be a $\IQ$-vector space of dimension $2n\geq4$
together with an endomorphism $\Phi$ and let
$H^*=\bigoplus_{k=0}^{4n} H^k$ be a graded $\IQ$-algebra. We
assume that there is a graded inclusion ${\bigwedge}^{2*}(A\oplus
A)\subset H^{2*}$ and  consider $B_1:=\bigwedge^2A\oplus\{0\}$,
$B_2:=\{0\}\oplus\bigwedge^2A$, and $H':=B_1\oplus B_2$
 as subspaces of $H^2$. We shall
use the notation of Corollary \ref{DelCor} with $\ell=2$,
$m=2n-1$.

\begin{hypo}\label{hypo3}
{\rm i)} $H^2=B_1\oplus B_2\oplus P$ with $P:=P_{2n-1}$ as in
(\ref{DFNP}),

{\rm ii)} $P=Q_1\oplus Q_2\oplus\IQ e_1\oplus \IQ e_2$ for some
vector spaces $Q_i$ and $P_1=0$,

{\rm iii)} $\{\alpha\in H^{2}_{\IC}~|~\alpha^2=0\}= \{\alpha\in
B_{1\IC}~|~\alpha^2=0\}\cup \{\alpha\in B_{2\IC}~|~\alpha^2=0\}$,

{\rm iv)} $\alpha^2a^{2n-2}=0$ for all $\alpha\in B_1$ and $a\in
P$, and

 {\rm v)} The kernel of the multiplication
$B_1\oplus B_2\stackrel{\cdot e_i}{\longrightarrow} H^4$,
$i=1,2$, is the diagonal $\Delta:=\{(a,a)\}\subset
\bigwedge^2A\oplus\bigwedge^2 A$ for $i=1$ and the graph
$\Gamma_{\bigwedge^2\Phi}$ for $i=2$. Similarly, ${\rm Ker}(\cdot
a_i)=B_i$ for any $0\ne a_i\in Q_i$. The sum $Q_1\cdot
B_2+Q_2\cdot B_1+{\rm Im}(\cdot e_1)+{\rm Im}(\cdot e_2)$ is
direct.
\end{hypo}

\begin{prop}\label{Propnotproj2}
Suppose $H^*$ and $\Phi$ meet the requirements of \ref {hypo3} and
\ref{hypo}, respectively. Then $H^*$ cannot be realized as the
rational cohomology ring $H^*(X,\IQ)$ of a projective manifold
$X$.
\end{prop}

\begin{proof} Suppose $X$ is a projective manifold whose rational
cohomology ring $H^*(X,\IQ)$ can be identified with $ H^*$.

Due to iii) and Lemma \ref{DeligneLem}, $B_1$, $B_2$, and hence
$H'$ are Hodge substructures of $H^2$. Thus, Corollary
\ref{DelCor} applies and shows that $P$ is a Hodge substructure.
(Note that $\bigwedge^2 A$ is spanned by vectors $\alpha$ with
$\alpha^2=0$.)

Due to v), the algebraic set $Z\subset P$ (see notation in
Corollary \ref{DelCor}) contains $\IC e_1$ and $\IC  e_2$
as two irreducible components. Indeed, if $\sum
a_i+\sum\eta_ie_i\in Z$ with $a_i\in Q_i$, then  some $0\ne
b=b_1+b_2\in B_1\oplus B_2$ is annihilated by it. Since the sum of
the multiplications is direct, this yields $a_2\cdot b_1=a_1\cdot
b_2=\eta_i(b\cdot e_i)=0$. In particular, $a_1\ne0$ implies
$b_2=0$ and $a_2\ne0$ implies $b_1=0$. Thus, if $\eta_1=\eta_2=0$,
then either $a_1\ne0$ or $a_2\ne0$.  Similarly, if
$\eta_1\ne0=\eta_2$, then $b_1=b_2\ne0$ and, therefore,
$a_1=a_2=0$. Finally, the case $\eta_1\ne0\ne\eta_2$ is excluded
by $\Delta\cap\Gamma_{\bigwedge^2\Phi}=\{(0,0)\}$, which follows
from $\mu_i\cdot \mu_j\ne1$ for all $i\ne j$ and $n\geq2$. (The
argument shows that the other irreducible components are
$Q_{i\IC}$.)

 Thus, by iii) of Corollary
\ref{DelCor}, the diagonal and the graph of $\bigwedge^2\Phi$ are
Hodge substructures of $B_1\oplus B_2$. In other words,
$\bigwedge^2\Phi$ is an endomorphism of the Hodge structure of
$\bigwedge^2 A$ induced by $B_1$ (or, equivalently, by $B_2$).

Clearly, $\bigwedge^2 A$ contains a subspace $V$ of dimension at
least two such that  $0=\alpha^2\in H^4$ for all $\alpha\in V$.
(For instance, take $V=\langle v_1\wedge v_2, v_1\wedge
v_3\rangle$ if $A=\bigoplus \IQ v_i$.)

Hence, by the Hodge--Riemann bilinear relations this excludes
$V\subset H^{1,1}(X)$ (see Example \ref{exaHSwttwo}). Therefore,
$\bigwedge^{2,0}A\ne0$ and, hence, the Hodge structure
$\bigwedge^2 A$ does not contain any Hodge class (see Proposition
\ref{PropnoHodgeclwttwo}).

This shows that all Hodge classes of $H^2$ are contained in $P$.
In particular, any hyperplane class $[\omega]$ is contained in
$P$. On the other hand, due to iv) one has
$\alpha^2.[\omega]^{2n-2}=0$ for all $\alpha\in B_1$, but
$H^2(X,\IQ)$ can clearly not contain a Hodge substructure of
dimension $\geq2$  which is isotropic with respect to the
polarization (see Example \ref{exaHSwttwo}). This yields the
contradiction.
\end{proof}
\section{Construction of examples}\label{Exassect}

So far we have explained how Voisin is able to exclude certain
Hodge structures on $\IQ$-algebras from being realized by the
cohomology of a projective manifold. It remains to find compact
K\"ahler manifolds which do realize these structures and which,
therefore, are topologically different from any projective
manifold.

The first two examples are obtained as blow-ups of well-known
K\"ahler manifolds and the following general facts will be used
tacitly throughout (see \cite{Demailly,GH,V4}). Let
$\pi:\widetilde X\to X$ be the blow-up of a compact complex
manifold $X$ along a submanifold $i:Y\hookrightarrow X$ of
codimension $c\geq2$. The exceptional divisor
$j:E=\pi^{-1}(Y)\hookrightarrow \widetilde X$ is isomorphic to
$\IP({\mathcal N}_{Y/X})$ and $\pi|_E$ equals the projection
$\pi_Y:\IP({\mathcal N}_{Y/X})\to Y$. In the following,
cohomology  will be considered with coefficients in $\IQ$.

$\bullet$ If $X$ is K\"ahler, then $\widetilde X$ is K\"ahler.

$\bullet$ If a submanifold $Z\subset X$ intersects $Y$
transversally, then the proper transform, which is by definition
the closure of $\pi^{-1}(Z\setminus Y)$, is the blow-up
$\widetilde Z\to Z$ along $Y\cap Z$.

 $\bullet$ The natural morphisms $\pi^*:H^k(X)\to
H^k(\widetilde X)$ and
$$H^{k-2(\ell+1)}(Y)\stackrel{\pi_Y^*}{\longrightarrow}H^{k-2(\ell+1)}(E)
\stackrel{\cdot
h^\ell}{\longrightarrow}H^{k-2}(E)\stackrel{j_*}{\longrightarrow}H^k(\widetilde
X),$$ where $h:={\rm c}_1({\mathcal O}_{\pi_Y}(1))$, induce
isomorphisms
$$H^k(\widetilde X)\cong
H^k(X)\oplus\bigoplus_{i=k-2(c-1)}^{k-2}H^i(Y).$$ In particular,
$H^2(\widetilde X)\cong H^2(X)\oplus\IQ e$ if $e:=[E]\in
H^2(\widetilde X)$ and $Y$ is connected.

$\bullet$ Moreover,
$$\varphi_e:H^k(X)\stackrel{\pi^*}{\longrightarrow}
H^k(\widetilde X)\stackrel{\cdot e}{\longrightarrow}
H^{k+2}(\widetilde X)$$ equals
$$H^k(X)\stackrel{i^*}{\longrightarrow} H^k(Y)\stackrel{\pi_Y^*}{\longrightarrow}
H^k(E)\stackrel{j_*}{\longrightarrow} H^{k+2}(\widetilde X).$$ In
particular, ${\rm Ker}(H^k(X)\stackrel{\varphi_e}{\longrightarrow}
H^{k+2}(\widetilde X))={\rm
Ker}(H^k(X)\stackrel{i^*}{\longrightarrow} H^k(Y))$.

$\bullet$ If $Y=Y_1\sqcup Y_2$ and accordingly $E=E_1\sqcup E_2$,
then for $k=1$ the sum $\sum{\rm Im}(\varphi_{e_i})\subset
H^3(\widetilde X)=H^3(X)\oplus H^1(Y_1)\oplus H^1(Y_2)$ is direct
and similar for $k=2$ the sum $\sum{\rm Im}(\varphi_{e_i})\subset
H^4(\widetilde X)\cong H^4(X)\oplus\bigoplus
H^2(Y_i)\oplus\bigoplus H^0(Y_i)$ is direct. (Note that the
degree zero terms only occur if $c\geq3$.) This principle can be
generalized to the case that $Y_1,Y_2$ intersect transversally
and that $\pi:\widetilde X\to X$ is obtained from first
blowing-up along $Y_1$ and then along the proper transform of
$Y_2$.
\subsection{Voisin's first example}\label{FirstEx}

Let $\Phi$ be an endomorphism of a $\IQ$-vector space $A$ of dimension
$2n\geq 4$
satisfying Hypothesis \ref{hypo}. By passing to $k\Phi$ for some $0\ne k\in\IZ$
if necessary, we may assume that $\Phi^*$ preserves a maximal lattice
$\Gamma\subset A^*$. Consider the complex torus $T:={A^{1,0}}^*/\Gamma$,
where $A_\IC=A^{1,0}\oplus
A^{0,1}$ is a Hodge structure as in Example \ref{TorusPhi}.
Then there exist natural isomorphisms $H^1(T,\IQ)\cong A$ and
$H^{1,0}(T)\cong A^{1,0}$.
The endomorphism $\Phi^*$ induces an endomorphism of $T$ which shall also
be denoted $\Phi^*$.

\begin{rema}
The complex tori $T$ and $T\times T$ are not projective due to Corollary
\ref{CorwtonenoPo}, but they are, as all other complex tori, deformation
equivalent and hence homeomorphic to abelian varieties.
\end {rema}

Voisin's first example constructed in \cite{V1} is a compact K\"ahler
manifold $X$ obtained as a blow-up of $T\times T$.

Consider the following submanifolds of $T\times T$:
$$\Delta_1:=\{(x,-x)\},~~\Delta_2:=\{(x,-\Phi^*(x))\},
~~T_1:=\{0\}\times T,~~T_2:=T\times\{0\},$$
which meet pairwise transversally. (E.g., via the first
projection the tangent space of $\Delta_1\cap\Delta_2$ in an
intersection point $z=(x,y)$ is identified with ${\rm Ker}({\rm
id}-\Phi^*)$, but $1$ is not an eigenvalue of $\Phi$.)

Let $z_1,\ldots,z_M\in T\times T$ be the finitely many
intersection points of all the pairwise intersections. Then
consider the blow-up $\pi_1:\widetilde{T\times T}\to T\times T$
in these points. The proper transforms of the four submanifolds
$\widetilde\Delta_1,\widetilde\Delta_2, \widetilde T_1,\widetilde
T_2$ are pairwise disjoint sub\-mani\-folds of
$\widetilde{T\times T}$. Thus, the blow-up
$\pi_2:X\to\widetilde{T\times T}$ along the union
$\widetilde\Delta_1\cup\widetilde\Delta_2\cup\widetilde
T_1\cup\widetilde T_2$ is a compact K\"ahler manifold. 

We shall denote by $F_1,\ldots, F_M\subset X$ the proper
transform of the exceptional divisors of $\pi_1$ and by
$E_1\to\widetilde T_1$, $E_2\to\widetilde T_2$,
$E_3\to\widetilde\Delta_1$, $E_4\to\widetilde\Delta_4$ the
exceptional divisors of $\pi_2$. Their cohomology classes shall
be called $f_1,\ldots,f_M,e_1,\ldots ,e_4\in H^2(X,\IQ)$. It is
the second blow-up $\pi_2$ and its exceptional classes
$e_1,\ldots,e_4$ that are important; the first blow-up $\pi_1$ is
only needed in order to ensure the smoothness of $X$.

The composition $\pi:=\pi_1\circ\pi_2:X\to T\times T$ induces a
graded inclusion $\bigwedge^*(A\oplus A)=H^*(T\times T,\IQ)\subset
H^*(X,\IQ)$.

\begin{prop}
The conditions {\rm i)-iv)} of  \ref{hypo2} are satisfied.
\end{prop}

\begin{proof}
The condition i) is obvious, as $X$ and $T\times T$ are
homeomorphic away from subsets of real codimension $\geq2$.
Since $H^2(T\times T,\IQ)\cong\bigwedge^2 H^1(T\times T,\IQ)$,
one has $H^2(X,\IQ)\cong\bigwedge^2(A\oplus
A)\oplus\bigoplus_{i=1}^M\IQ f_i\oplus\bigoplus_{i=1}^4\IQ e_i$.

A class in $\bigwedge^{4n-2}
H^1(X,\IQ)=\bigwedge^{4n-2}H^1(T\times T,\IQ)$ can be thought of
as a linear combination of fundamental classes of subsets of real
codimension $4n-2$ in $T\times T$ in general position, whose
pull-back clearly avoids the exceptional divisors $F_1,\ldots,
F_M,E_1,\ldots,E_4$ which all live over subsets of real
codimension $>2$. This yields ii) with $P=\langle
f_1,\ldots,f_M,e_1,\ldots,e_4\rangle$ and $R=0$.

A similar argument yields iii), where $P_1=\langle f_1,\ldots,
f_M\rangle$. Finally, condition iv) is proved by applying the
above general remarks on the cohomology of a blow-up and by using
the explicit description of $\Delta_1,\Delta_2,T_1$, and $T_2$.
\end{proof}

Together with Proposition \ref{Propnotproj}  this yields
\begin{coro}
The rational homotopy type of the compact K\"ahler manifold $X$
of dimension $2n\geq4$ is not
realized by any projective manifold.\qqed
\end{coro}
Note that this time the result has been phrased in terms of the
rational homotopy type rather than in terms of the rational
cohomology. Both statements are equivalent due to \cite{DGMS} and
the fact that the fundamental group is abelian in our situation.

\begin{rema}\label{avoidP'Rem}
One could also avoid the initial point blow-ups and instead
successively blow-up $T_1$, $T_2$, $\Delta_1$, $\Delta_2$,
respectively their proper transforms. The above arguments remain
valid, only that in this case $P_1=0$.
\end{rema}

In order to fully prove Theorem \ref{VoisinMainThm} it remains to
construct examples of odd dimension. These are obtained as
products $X':=X\times \IP^1$, where $X$ is one of the compact
K\"ahler manifolds above. Once more the conditions {\rm i)-iv)} of
\ref{hypo2} are satisfied, but this time $R=H^2(\IP^1,\IQ)$. The
rest of the argument is unaffected by this modification.

\begin{rema}
In \cite{V1} it is first shown that the integral cohomology
$H^*(X,\IZ)$ of the above constructed K\"ahler manifold cannot be
realized by a projective manifold. The proof of this weaker
statement does not rely on Deligne's principle, but uses the
Albanese morphism instead.

One finds in \cite{V1} also an example, due to Deligne, of a
compact K\"ahler manifold whose complex cohomology $H^*(X,\IC)$
cannot be realized by a projective manifold. The manifold $X$ is
again obtained as a blow-up of $T\times T$.
\end{rema}

\subsection{Simply-connected examples}

One might wonder whether the fundamental group is responsible for
the fact that the above constructed compact K\"ahler manifold is
topologically different from any projective manifold. This
question lead Voisin to her second example, which is
simply-connected. Roughly, the simply-connected K\"ahler manifold
is obtained from the first one by dividing by the
$\IZ/2\IZ\times\IZ/2\IZ$-action, which is induced by the standard
involution on the two factors.

On the one hand, the construction is simpler in the sense that
blowing-up $T_1$ and $T_2$ can be avoided, which was needed before
to detect certain Hodge substructures. As it turns out, the
analogous Hodge structures in the simply-connected case can be
described directly. (As the examples will be simply-connected,
one cannot work with Hodge structures of weight one. Therefore,
Voisin analyses the weight-two Hodge structure on $H^2(X,\IQ)$
instead.) On the other hand, due to the (mild) singularities of
$T/\pm$, the construction is slightly more involved, as we first
have to desingularize.

In \cite{V1} Voisin proceeds as follows. Start with a torus
$T={A^{1,0}}^*/\Gamma$ as in Section \ref{FirstEx}. In
particular, $T$ comes with an endomorphism $\Phi^*$. Next,
consider the quotient $T/\pm$ of $T$ by the standard involution
$z\mapsto \pm z$ and its desingularization $K\to T/\pm$ obtained
by a simple blow-up of all the two-torsion points. Equivalently,
one may first blow-up the two-torsion points  $\widetilde T\to T$
and then take the quotient $K=\widetilde T/\pm$ by the induced
involution. The latter description shows that $K$ is smooth and
K\"ahler. (Indeed, a general result of Varouchas \cite{Varou}
proves that for a surjection $\pi:X\to X'$ whose fibres are all
of dimension $\dim(X)-\dim(X')$ the manifold $X'$ is K\"ahler if
$X$ is so.) Viewing $K$ as the desingularization of $T/\pm$,
shows that it is simply-connected, for $T/\pm$ is.

The endomorphism $-\Phi^*$ of $T$ descends to an endomorphism
$-\bar\Phi^*$ of $T/\pm$ and we consider its graph
$\Gamma_{-\bar\Phi^*}\subset (T/\pm)\times (T/\pm)$.

In the last step, one first blows-up $K\times K$ along the
anti-diagonal $\Delta_1:=\{(a,-a)\}$ and then along the proper
transform $\Gamma'$ of $\Gamma_{-\bar\Phi^*}$. (Note that
$\Gamma'$ is smooth. This can be seen by passing via $\widetilde
T\times\widetilde T\to T\times T$.)

 Thus, the resulting variety $X$ is indeed a
K\"ahler manifold. We let $\pi:X\to K\times K$ be the composition
of the two blow-ups. The two exceptional divisors $E_1\to \Delta$
and $E_2\to\Gamma'$ yield distinguished cohomology classes
$e_1,e_2\in H^2(X,\IZ)$.

\begin{prop} Let $n\geq3$.
Then the conditions {\rm i)-v)} of \ref{hypo3} are satisfied.
\end{prop}

\begin{proof}
Since the involution of $T$ acts trivially on $H^2(T,\IQ)$,
one has $H^2(T/\pm,\IQ)\cong H^2(T,\IQ)=A$ and
$H^2(K,\IQ)=A\oplus\bigoplus \IQ f_j$, where $f_i$ are the
classes corresponding to the exceptional divisors $F_i$ over
the two-torsion points.

Thus, $H^2(X,\IQ)=H^2(K\times K,\IQ)\oplus\IQ e_1\oplus \IQ e_2
=H^2((T/\pm)\times (T/\pm),\IQ)\oplus Q_1 \oplus Q_2 \oplus\IQ
e_1\oplus \IQ e_2 $, where $Q_i$ is the pull-back of
$\bigoplus\IQ f_j$ under the $i$-th projection onto $K$.

It is easy to see that $P:=Q_1\oplus Q_2\oplus\IQ e_1\oplus\IQ
e_2$ is indeed the subspace that is annihilated by
$S^{2n-1}H^2((T/\pm)\times (T/\pm),\IQ)$. This proves i).

Since $\bigwedge^2 A$ is spanned by elements $a$ with $a^2=0$ and
no non-trivial linear combination of
$f_{1j}:=\pi_1^*f_j,f_{2j}:=\pi_2^*f_j,e_1$, and $e_2$ has this
property, condition iii) follows. It is here that one needs the
assumption $n\geq3$. The verification of condition v) is
straightforward; use the explicit description of the classes $e_1$
and $e_2$.

To conclude, we have to verify condition iv). One can show that
for all $\alpha\in B_1$ expressions of the form $\alpha^2\cdot
P(f_{ij},e_1,e_2)$ with $P$ a polynomial of degree $2n-2$ are
indeed trivial. Here are a few of the necessary arguments.
Firstly, $f_{ij}^k=0$ for all $k>n$. Secondly, the classes
$f_{ij}\cdot e_k$ and $e_1\cdot e_2$ are supported over finitely
many points in $(T/\pm)\times (T/\pm)$ and, hence as $\alpha$ is
pulled-back from there, one has $\alpha \cdot(f_{ij}\cdot
e_k)=\alpha\cdot (e_1\cdot e_2)=0$. Thirdly, $\alpha\cdot
f_{1j}=0$. Thus, the only combinations that need to be checked
are $\alpha^2\cdot e_i^{2n-2}$. We may assume that
$E_i=\IP(\Omega_T)$ and that $\pi|_{E_i}$ is the natural
projection $p:E\to T$. Then one shows that $e_i|_{E_i}={\rm
c}_1({\mathcal O}_p(-1))$ and thus reduces to
$0=p^*\alpha_T^2.{\rm c}_1({\mathcal O}_p(-1))^{2n-3}$, which
follows from ${\rm c}_1({\mathcal O}_p(-1))^{k}=0$ for $k\geq n$
and the assumption $n\geq3$.
%
\end{proof}

Together with Proposition \ref{Propnotproj2} this yields

\begin{coro}
The rational homotopy type of the compact simply-connected K\"ahler
manifold $X$ of dimension $2n\geq6$
is not realized by any projective manifold.\qqed
\end{coro}

Odd-dimensional examples can again be produced by taking products with
$\IP^1$. In \ref {hypo3} only i) and iii) have to be modified. In
i) one has $H^2=B_1\oplus B_2\oplus P\oplus R$ with $R=H^2(\IP^1,\IQ)$
and in iii) $R_\IC$ will provide another irreducible component. The arguments are not
affected by this modification. This yields
  C.\ Voisin's second counter-example:
\begin{theo}[\cite{V1}] In any dimension $\geq6$ there exists
a simply-connected compact K\"ahler manifold which does not have
the rational homotopy type of a projective manifold.
\end{theo}

Once more, instead of working with the rational homotopy type one
could equivalently say that $H^*(X,\IQ)$ is not realized as the
cohomology ring of a projective manifold (see \cite{DGMS}).

\begin{rema}
Inspired by Voisin's examples, Oguiso studies in \cite{Oguiso}
simply-connected compact K\"ahler manifolds of dimension $d\geq4$
which are not projective, but rigid, i.e.\ which do not allow any
deformations at all and, therefore, cannot be deformed to
projective ones in particular. In the case of simply-connected
examples one can no longer work with Hodge structures of weight
one. Thus, K3 surfaces (or, more generally, compact hyperk\"ahler
manifolds) with their very special but rich Hodge structures of
weight two provide a  reservoir of potentially interesting
examples. Roughly, the special endomorphisms of tori used by
Voisin are in \cite{Oguiso} replaced by special automorphisms of
K3 surfaces which are described completely by their action on the
second cohomology.

However, the methods in  \cite{Oguiso} fall short of proving that
the examples do not have the rational homotopy type of projective
manifolds. It seems likely, nevertheless, that four-dimensional
simply-connected examples could eventually be produced in this
way.
\end{rema}

\subsection{The birational Kodaira problem}

Right after \cite{V1} had appeared,  modifications of the original
problem have been proposed.  For many problems in complex
algebraic geometry it is natural not to restrict to projective or
K\"ahler manifolds, but to allow manifolds that are birational or
bimeromorphic to those. Passing to a bimeromorphic model often
changes the topology drastically, but in a somewhat controlled
manner. So, modifying Kodaira's problem in this sense seems
natural also from a topological point of view.

More precisely, the compact K\"ahler manifolds constructed in
\cite{V1} are both bimeromorphic to compact K\"ahler manifolds
which do have the homotopy type of projective manifolds. E.g.\ in
the first example, described in Section \ref{FirstEx}, the
K\"ahler manifold $X$ was constructed as a blow-up of a torus
whose underlying manifold carries also the structure of a
projective manifold. In other words, after a controlled
topological modi\-fication the original topological manifold
underlying $X$ has been transformed to one that does carry a
projective structure. So, one could ask whether this is true for
any K\"ahler manifold. Again, the answer is negative.

\begin{theo}[\cite{V2}] There exist compact K\"ahler manifolds
$X$ of dimension $2n\geq10$ such that no complex manifold
bimeromorphic to it has the rational homotopy type of a
projective manifold.
\end{theo}

The principal ideas in \cite{V2} are similar to those in
\cite{V1}. Roughly, one tries to detect certain Hodge structures
in terms of the multiplicative structure of the cohomology ring
and to derive a contradiction to the existence of a polarization
on the (primitive) second cohomology of a projective manifold.
Technically, the arguments are more involved and we only give an
idea of the actual construction.

The construction of the birational counter-examples in \cite{V2}
starts again with the same torus $T$ of dimension $n\geq 4$ and an
endomorphisms $\Phi$ satisfying \ref{hypo}. If ${\mathcal P}$
denotes the Poincar\'e bundle on $T\times\widehat T$, then let
$E:={\mathcal P}\oplus{\mathcal P}^{-1}$ and $E_\Phi:=(\Phi,{\rm
id})^*E$. In the next step one considers the fibre product
$\IP(E)\times_{T\times\widehat T}\IP(E_\Phi)$ and its quotient
$Q$ by the action of $(\IZ/2\IZ)\times(\IZ/2\IZ)$ given by natural
lifts of $(-{\rm id},{\rm id})$ and $({\rm id},-{\rm id})$. Then
any K\"ahler desingularization $X$ of $Q$ will work. Note that
these examples are bimeromorphic to a $\IP^1\times\IP^1$-bundle
over $K\times \widehat K$, where $K\to T/\pm$ is the
desingularization considered in the simply-connected case.

The reason that one is able to control in this example all
bimeromorphic models by cohomological methods is due to the fact
that there exist only few subvarieties of positive dimension.

\smallskip

\section{Further comments}

This is still not the end. Why not allowing topological changes
that are not obtained by bimeromorphic maps? One could ask whether
there always exists another complex structure on $X$ (e.g.\ one
obtained by a deformation) such that a bimeromorphic model of
this new one has the rational homotopy type of a projective
manifold. So, more formally, if one introduces the equivalence
relation between complex manifolds generated by deformations and
bimeromorphic correspondences, one might ask whether any compact
K\"ahler manifold is equivalent to a projective manifold.

Continuing in this direction, one could allow singular varieties
or certain ramified covers in order to enlarge the equivalence
classes. Would the answer to Kodaira's problem be different
then?  Most of these questions are open for the time being, but
see the comments in \cite{V3}.

In another direction, it could be interesting to see whether
the birational geometry does matter in these questions. The above
counter-example for the birational Kodaira problem is, by construction,
of Kodaira dimension $-\infty$. For the time being the techniques do
not seem to produce examples of non-negative Kodaira dimension.

As has been mentioned, topologically there is no difference
between compact K\"ahler surfaces and projective surfaces. Due to
the examples of Voisin, the situation changes drastically in
dimension $\geq4$ (or rather $\geq6$ if one prefers
simply-connected manifolds). What seems open, however, is the
three-dimensional case:

{\it Does there exist a compact K\"ahler threefold which
is not homeomorphic to a projective manifold?}

\smallskip
Since we mentioned fundamental groups in the beginning, let us
point out that the following problem is also still open:

{\it Does there exist a group that is the fundamental group of a
compact K\"ahler manifold, but not of a projective manifold?}

\smallskip

A question of a more general nature is the following:

{\it Are there topological, cohomological,... conditions that
decide whether a compact K\"ahler manifold can also be endowed
with a complex structure which is projective?}

Nothing seems to be known in this direction and the examples show
that if such conditions can be found at all, they cannot be
formulated purely in terms of the fundamental group.

\end{document}